\documentclass[10pt]{amsart}
\usepackage[cp1251]{inputenc}
\usepackage[english,russian]{babel}
\usepackage{amsmath}
\usepackage{amssymb}
\usepackage{amsfonts}

\def\udcs{512.5} 
\def\mscs{17A36} 
\newtheorem{lemma}{Лемма}
\newtheorem{theorem}{Теорема}

\newtheorem{corollary}{Следствие}
\newtheorem{proposition}{Предложение}

\def\logo{{\bf\huge S\raisebox{0.2ex}{\hspace{0.55ex}\raisebox{0.05ex}e\hspace{-1.65ex}$\bigcirc$}MR}}

\def\semrtop
     {
  \vbox{
     \noindent\logo
     \hspace{80mm}\raisebox{1ex}{ISSN 1813-3304 }

     \vspace{5mm}

     \begin{center}
     {\huge СИБИРСКИЕ \ ЭЛЕКТРОННЫЕ} \\[2mm]
     {\huge МАТЕМАТИЧЕСКИЕ ИЗВЕСТИЯ} \\[2mm]
     {\large Siberian Electronic Mathematical Reports} \\[1mm]
     {\LARGE\tt{http://semr.math.nsc.ru}}\\[0.5mm]
     \end{center}
     \vspace{-3mm}
     \noindent
     \begin{tabular}{c}
     \hphantom{aaaaaaaaaaaaaaaaaaaaaaaaaaaaaaaaaaaaaaaaaaaaaaaaaaaaaaaaaaaaaaaaaaaaaa} \\
     \hline\hline
     \end{tabular}

     \vspace{1mm}
     {\flushleft\it Том 15, стр. 144--153 (2017) \hspace{65mm}{\rm\small УДК \udcs}}
     \newline
     	{\rm\small DOI~10.17377/semi.2017.12.xxx}\hphantom{aaaaaaaaaaaaaaaaaaaaaaaaaaaaaaaaaaaa}{\rm\small MSC\ \ \mscs }
  }
}

\begin{document}
\renewcommand{\refname}{References}

\thispagestyle{empty}

\title[Автоморфизм Нагаты свободных неассоциативных алгебр]{Автоморфизм Нагаты свободных неассоциативных алгебр ранга два над евклидовыми кольцами}

\author{{А.А. АЛИМБАЕВ}, {У.У. УМИРБАЕВ}}%
\address{Alibek Alimbaev  
\newline\hphantom{iii} Kostanay State Pedagogical Institute,
\newline\hphantom{iii} Tauelsizdik street, 118,
\newline\hphantom{iii} 110000, Kostanay, Kazakhstan}%
\email{alialimbaev@gmail.com}%

\address{Ualbai Umirbaev  
\newline\hphantom{iii} Wayne State University,
\newline\hphantom{iii} 656 W. Kirby,
\newline\hphantom{iii} Detroit, MI 48202, USA}%
\email{umirbaev@math.wayne.edu}%

\thanks{\sc Alimbaev A.A., Umirbaev U.U,
The Nagata automorphism of free nonassociative algebras of rank two over Euclidean domains}
\thanks{\copyright \ 2017 Алимбаев А.А., Умирбаев У.У}
\thanks{\rm Работа поддержана МОН РК (грант 0538/GF4)}
\thanks{\it Поступила 3 октября 2017 г., опубликована 31 декабря 2017 г.}%

\semrtop \vspace{1cm}
\maketitle {\small
\begin{quote}
\noindent{\sc Abstract.} We construct an analogue of the Nagata automorphism of free nonassociative algebras and free commutative algebras of rank two over а Euclidean domain and prove that it is wild.\medskip

 \end{quote}
}

\section{Введение}
\hspace*{\parindent}

В 1942 году Х.В.Е. Юнг \cite{1} доказал, что все автоморфизмы алгебры многочленов ${\bf k}[x,y]$ от двух переменных над полем {\bf k} характеристики 0 являются ручными. В 1953 году В. ван дер Калк \cite{2} обобщил этот результат для полей произвольной характеристики.

В 1972 году М. Нагата \cite{3} построил автоморфизм
$$
\sigma=(x+2yw+w^2z,y+wz,z), w=xz-y^2
$$
алгебры многочленов ${\bf k}[x,y,z]$ и высказал гипотезу, что этот автоморфизм не является ручным (т.е. является диким). В 2004 году У. Умирбаев и И. Шестаков доказали \cite{4, 5, 6}, что автоморфизм Нагаты является диким в случае полей нулевой характеристики.

В работах Л. Макар-Лиманова \cite{7} и А. Чернякевич \cite{8} доказано, что автоморфизмы свободных ассоциативных алгебр ранга 2 являются ручными. Кроме того, ими было доказано, что группы автоморфизмов алгебр многочленов ${\bf k}[x,y]$ и свободной ассоциативной алгебры ${\bf k}\left\langle x,y\right\rangle$ от двух переменных изоморфны.

Автоморфизмы двупорожденных правосимметричных алгебр над произвольными полями \cite{9} и двупорожденных свободных алгебр Пуассона  над полями нулевой характеристики  \cite{10} также являются ручными.

В случае свободных ассоциативных алгебр ранга три вопрос о существовании диких автоморфизмов (проблема Кона) был решен У.У. Умирбаевым. Им было доказано \cite{11}, что автоморфизм Аника
$$
\delta=(x+z(xz-zy), y+(xz-zy)z, z)
$$
свободной ассоциативной алгебры ${\bf k}\left\langle x,y,z\right\rangle$ над полем характеристики 0 является диким.

В 1964 году П. Кон \cite{12} доказал, что автоморфизмы конечно порожденных свободных алгебр Ли над произвольными полями являются ручными.
В 1968 году Дж. Левин \cite{13} обобщил этот результат для шрайеровых многообразий алгебр. Напомним, что шрайеровыми являются многообразия всех неассоциативных алгебр \cite{14}, коммутативных и антикоммутативных алгебр \cite{15}, алгебр Ли \cite{16, 17} и супералгебр Ли \cite{18,19}.
Следовательно, автоморфизмы свободных неассоциативных алгебр, свободных коммутативных и антикоммутативных алгебр конечного ранга над полями также являются ручными.

Группы автоморфизмов конечно порожденных свободных алгебр над кольцами главных идеалов были исследованы в работах \cite{20}, \cite{21}. В частности, один из результатов (Теорема 3, \cite{21}) гласит, что автоморфизмы свободных неассоциативных алгебр над кольцами главных идеалов являются ручными. 

В данной работе нами построен пример дикого автоморфизма свободной неассоциативной алгебры и коммутативной алгебры ранга два над евклидовыми кольцами. Этот пример является аналогом автоморфизма Нагаты.
Статья организована следующим образом. В разделе 2 мы приведем необходимые определения и известное представление группы ручных автоморфизмов алгебры $\Phi[x_1, x_2]$ над произвольной областью целостности $\Phi$ в виде свободного произведения подгрупп с объединенной подгруппой и сформулируем некоторые полезные утверждения о степени ручных автоморфизмов. В разделе 3 докажем, что всякий ручной автоморфизм алгебры $\Phi[x_1, x_2]$ над евклидовым кольцом $\Phi$ является элементарно сократимым. В  разделе 4 построен дикий автоморфизм свободных неассоциативных алгебр и свободных коммутативных алгебр ранга 2 над евклидовым кольцом $\Phi$.

\section{Представление автоморфизмов алгебры $\Phi[x_1, x_2]$}
\hspace*{\parindent}
Пусть $\Phi$ --- произвольная область целостности и $P=\Phi[x_1,x_2]$ --- алгебра многочленов от двух переменных над $\Phi$.  Множество всех обратимых элементов $\Phi$ обозначим через $\Phi^*$.

Пусть $Aut(P)$ --- группа всех автоморфизмов алгебры $P$. Обозначим через $
{\phi=(f_1,f_2)}
$ автоморфизм алгебры $P$ такой, что $\phi(x_i)=f_i,\, 1\leq i\leq 2$.

 Автоморфизмы вида
\begin{flushleft}
$$
\delta_1:\begin{cases}
x_1 \mapsto \alpha x_1+ f(x_2)\\
x_2 \mapsto x_2\\
\end{cases},\ \
\delta_2:\begin{cases}
x_1 \mapsto x_1\\
x_2 \mapsto \beta x_2+ g(x_1)\\
\end{cases},
$$
\end{flushleft}
где $\alpha, \beta \in \Phi^*$ называются \textit{элементарными}. Подгруппа $T(P)$ группы $Aut(P)$, порожденная всеми элементарными автоморфизмами, называется \textit{подгруппой ручных автоморфизмов}.
Автоморфизм $\varphi\in Aut(P)$ называется {\em ручным}, если $\varphi\in T(P)$, иначе $\varphi $ называется {\em диким.}

{\em Элементарным преобразованием} системы элементов $(f_1,f_2)$ называется замена одного элемента $f_i$ на элемент вида $\alpha f_i +g,$ где $\alpha \in \Phi^*, g\in \left\langle {f_j |j\ne i}\right\rangle.$

Запись
$$
(f_1, f_2)\to (g_1, g_2)
$$
означает, что система элементов $(g_1, g_2)$ получена из системы элементов $(f_1, f_2)$ одним элементарным преобразованием.

Если $(f_1,f_2)$ --- ручной автоморфизм алгебры $P$, то существует последовательность элементарных преобразований вида

$$
(x_1, x_2)=(f^{(0)}_1,f^{(0)}_2)\rightarrow (f^{(1)}_1,f^{(1)}_2)\rightarrow \ldots \rightarrow (f^{(k)}_1, f^{(k)}_2)=(f_1,f_2).
$$

Если
$$
\varphi=(f_1, f_2), \: \psi=(g_1, g_2),
$$
то произведение в $Aut(P)$ определяется по формуле
$$
\varphi \circ \psi=(g_1(f_1,f_2), g_2(f_1,f_2)).
$$

Автоморфизм $\lambda$ алгебры $P$ называется {\em аффинным}, если
$$
\lambda=(a_1x_1+b_1x_2+c_1, a_2x_1+b_2x_2+c_2),
$$
\[\begin{vmatrix}
a_1 & b_1 \\
a_2 & b_2
\end{vmatrix} \in \Phi^*, \ a_i, b_i, c_i\in \Phi.\;\;
\]
Группу аффинных автоморфизмов алгебры $P$ обозначим через $Af_2(\Phi).$

Автоморфизм $\mu$ алгебры $P$ называется {\emтреугольным}, если
$$
\mu=(ax_1+h(x_2), bx_2+b_1), \;\;
$$
где $a,b \in \Phi^*, \ b_1 \in \Phi, \ h(x_2)\in \Phi[x_2].$
Через $Tr_2(\Phi)$ обозначим группу треугольных автоморфизмов алгебры $P=\Phi[x_1,x_2].$

Пусть {\bf k} --- произвольное поле.
Тогда $T({\bf k}[x_1,x_2])=Aut({\bf k}[x_1,x_2])$\cite{1, 2}.
Более того, $Aut({\bf k}[x_1,x_2])$ имеет хорошо известное представление в виде свободного произведения подгрупп с объединенной подгруппой:

\begin{theorem}
\cite{3} Группа автоморфизмов алгебры ${\bf k}[x_1,x_2]$ является свободным произведением подгрупп аффинных автоморфизмов $Af_2({\bf k})$ и треугольных автоморфизмов $Tr_2({\bf k})$ с объединенной подгруппой $H({\bf k})=Af_2({\bf k}) \cap Tr_2({\bf k}),$ т.е.
$$
Aut({\bf k}[x_1,x_2])=Af_2({\bf k})*_{H({\bf k})} Tr_2({\bf k}).
$$
\end{theorem}

\begin{corollary}
\cite{22} Пусть $\Phi$ --- произвольная область целостности. Группа ручных автоморфизмов $T(P)$ алгебры $P=\Phi[x_1,x_2]$
является свободным произведением подгрупп аффинных автоморфизмов $Af_2(\Phi)$ и треугольных автоморфизмов $Tr_2(\Phi)$ с объединенной подгруппой $H(\Phi) = Af_2(\Phi)\cap Tr_2(\Phi)$, т.е.
$$
T(P)=Af_2(\Phi)*_{H(\Phi)} Tr_2(\Phi).
$$
\end{corollary}

Приведем некоторые необходимые для нас факты из доказательств этих утверждений \cite{22, 23}.

Множество
$$
B_0=\left\{\tau=(x_1+h(x_2), x_2)| h(x_2)\in x_2^2 \Phi[x_2]\right\}
$$
 является системой представителей левых смежных классов группы $Tr_2(\Phi)$ по подгруппе ${Af_2(\Phi) \cap Tr_2(\Phi)}.
$ Заметим, что $B_0$ включает единицу при $h=0.$

Пусть $A_0$ --- произвольная фиксированная система представителей (включающая единицу) левых смежных классов группы $Af_2(\Phi)$ по подгруппе ${Af_2(\Phi) \cap Tr_2(\Phi)}.$

Тогда любой ручной автоморфизм $\phi \in T(P)$ однозначно представим
в виде
\begin{equation}\label{1}
\phi= \sigma_1 \circ \tau_1 \circ \sigma_2 \circ \tau_2 \circ \cdots \circ \sigma_k \circ \tau_k \circ \lambda,
\end{equation}
где $\sigma_i \in A_0, \sigma_2, \ldots, \sigma_k \neq id, \tau_i \in B_0, \tau_1, \ldots, \tau_{k}\neq id, \lambda \in Af_2(\Phi).$

Если $\phi=(f,g),$
то положим
$$\deg(\phi)=\max \{\deg(f), \deg(g)\}.$$

\begin{proposition}
\cite{22, 23} Пусть
$$
\psi = \sigma_1 \circ \tau_1 \circ \sigma_2 \circ \tau_2 \circ \cdots \circ \sigma_k \circ \tau_k = (f,g),
$$
где $\sigma_i \in A_0, \sigma_2, \ldots, \sigma_k \neq id, \tau_i \in B_0, \tau_1, \ldots, \tau_{k}\neq id,$
$
\deg(\tau_i)=n_i>1.
$\\
Тогда
$$
\deg(f)=n_1\cdots n_k,\  \deg(g)=n_1\cdots n_{k-1}, \  \deg(\psi)=n_1\cdots n_k.
$$
\end{proposition}

\section{О сократимости ручных автоморфизмов}
\hspace*{\parindent}
Обозначим через $\mathbb{Z}_+=\left\{0, 1, 2, \ldots, n, \ldots\right\}$ множество неотрицательных целых чисел.
 Целостное кольцо $\Phi$, не являющееся полем, называется {\em евклидовым} \cite{24}, если существует функция
$ \left|\cdot\right|:\Phi\backslash \left\{0\right\} \longrightarrow \mathbb{Z}_{+}$
(называемая нормой), удовлетворяющая следующим условиям:

E1) для любых $a,b\in \Phi\backslash \left\{0\right\}$, $\left|ab\right|\geq \left|a\right|$, причем равенство имеет место только тогда, когда элемент $b$ обратим;

E2) для любых $a,b\in \Phi$, где $b \neq 0$, существуют такие $q, r\in \Phi$, что $a=bq + r$ и либо $r=0$, либо $\left|r\right| < \left|b\right|.$

Положим $e=|1| \in \mathbb{Z}_+,$ где $1 \in \Phi$---единица кольца $\Phi.$
Имеем $|a|=e$ тогда и только тогда, когда $a\in \Phi^*.$

Основными примерами евклидовых колец являются кольцо $\mathbb{Z}$ целых чисел с абсолютным значением целых чисел и кольцо ${\bf k}[x]$ многочленов над полем ${\bf k}$ со степенью многочленов. Следовательно, $e=1$ в кольце $\mathbb{Z}$ и $e=0$ в кольце ${\bf k}[x].$

Пусть $\Phi[x_1,x_2]$ --- алгебра многочленов от двух переменных над евклидовым кольцом $\Phi.$
Введем линейный порядок $\geq$ на множестве базисных слов
$
{W=\left\{x_1^l x_2^m|l,m \in \mathbb{Z}_+ \right\}.}
$ Положим $x_1>x_2.$
Пусть $\deg$ --- обычная функция степени в $\Phi[x_1,x_2]$ и $\deg_{x_i}$ --- функция степени по $x_i.$
Если $u, v \in W,$ то положим  $u<v$, если выполняется одно из следующих условий:

(i) $\deg(u)<\deg(v);$

(ii) $\deg(u)=\deg(v),$ $\deg_{x_1}(u)<\deg_{x_1}(v)$.

Каждый ненулевой элемент $f\in \Phi[x_1,x_2]$ записывается однозначно в виде

$$
f=\alpha_1 w_1 +\alpha_2 w_2+\ldots+ \alpha_n w_n,\; 0\neq\alpha_i\in\Phi,  w_1>w_2>\ldots >w_n\in W.
$$

Слово $w_1$ называется старшим словом (мономом) $f$, а $\alpha_1$ называется старшим коэффициентом $f$. Обозначим их через $\overline{f}$ и $lc(f)$ соответственно. Через $f'=\alpha_1 w_1$ обозначим старший член элемента $f.$

Пусть $f$ --- произвольный элемент из $\Phi[x_1,x_2]$. Положим
$$
D(f)=(\overline{f},|lc(f)|)
$$
 и назовем $D(f)$ показателем элемента $f.$

Пусть $\varphi =(f_1,f_2)$ --- автоморфизм алгебры $\Phi[x_1,x_2]$.
Тогда показателем  автоморфизма $\varphi$ назовем $$D(\varphi)=(u, v, |lc(f_1)|+|lc(f_2)|)\in W^2\times Z_{+},$$
где
$\left\{u, v\right\}= \left\{\overline{f_1},\overline{f_2}\right\}$  и $u\geq v.$

Заметим, что $D$ --- инвариантно относительно перестановки компонент автоморфизма, т.е.
$$
D(f_1, f_2)=D(f_2, f_1).
$$

Имеем
$$
D(id)=(x_1, x_2, 2e),
$$
где $id=(x_1, x_2)$ --- тождественный автоморфизм.
Более того, $D(f_1, f_2)=(x_1, x_2, 2e)$ тогда и только тогда, когда элементы $f_1, f_2$ имеют следующий вид:
$$
f_1=\alpha_1x_1+ \beta_1x_2 + \gamma_1, \ f_2= \beta_2x_2+\gamma_2,
\eqno(2)
$$
или
$$
f_1= \beta_3x_2+\gamma_3, \ f_2=\alpha_2x_1+ \beta_4x_2 + \gamma_4,
\eqno(3)
$$
где $\alpha_1, \alpha_2,\beta_2, \beta_3\in\Phi^*$ и $\beta_1, \beta_4, \gamma_1, \gamma_2, \gamma_3, \gamma_4 \in\Phi.$

Обозначим через $\preceq$ лексикографический порядок на множестве $W^2\times Z_{+}$. Заметим, что множество $W^2\times Z_{+}$ линейно упорядочено относительно $\preceq$.

Автоморфизм $\theta=(f_1, f_2)$ называется {\em элементарно $D$-сократимым} (или $\theta$ {\em допускает элементарное $D$-сокращение}), если существует автоморфизм $\psi$ такой, что ${\theta \to \psi}$ и $D(\psi)\prec D(\theta).$
Будем говорить, что {\em автоморфизм $\psi$ является элементарным $D$-сокращением автоморфизма $\theta.$}
\begin{corollary}
Если автоморфизм $\theta=(f_1,f_2)$ допускает элементарное $D$-сокращение, то $f'_1 \in \left\langle f'_2\right\rangle$ или $f'_2 \in \left\langle f'_1\right\rangle.$
\end{corollary}
\begin{proof}
Допустим, что автоморфизм $\psi=(g_1, g_2)$ такой, что $\theta \to \psi$ и ${D(\psi)\prec D(\theta),}$ элементарно $D$-сокращает элемент $f_1$ автоморфизма $\theta.$
Следовательно,
$$
g_2 = f_2, \: g_1= \alpha f_1 +s(f_2),
$$
где $\alpha \in \Phi^*, s(f_2) \in \left\langle f_2 \right\rangle.$
Так как $D(g_1)\prec D(f_1),$ то $\alpha f'_1+s(f_2)'=0.$ Следовательно,
$\alpha f'_1+s'(f'_2)=0 \Rightarrow -\alpha f'_1=s'(f'_2) \Rightarrow f'_1 \in \left\langle f'_2\right\rangle.$ Аналогично, если автоморфизм $\psi$ элементарно $D$-сокращает элемент $f_2$ автоморфизма $\theta,$ то $f'_2 \in \left\langle f'_1\right\rangle.$
\end{proof}

\begin{lemma}\label{sl}
Пусть $\pi=(h_1, h_2)$ --- автоморфизм алгебры $\Phi[x_1,x_2].$ Если  старшие слова элементов $h_1, h_2$ равны, то автоморфизм $\pi$ является элементарно $D$-сократимым.
\end{lemma}
\begin{proof}
Так как $\overline{h_1} =\overline{h_2},$ то без ограничения общности предположим, что $|lc(h_1)|\geq |lc(h_2)|.$
По E1) существуют  $q,r\in \Phi$ такие, что ${lc(h_1)= lc(h_2) q+r}$ и либо $r=0,$ либо $|r|<|lc(h_2)|.$
Рассмотрим элементарное преобразование
$$\pi=(h_1, h_2)\to (h_1-q h_2, h_2)=\delta.$$
Имеем $D(h_1)\succ D(h_1-q h_2).$
Следовательно, $D(\pi)\succ D(\delta)$ и автоморфизм $\pi$ является элементарно $D$-сократимым.
\end{proof}

Характеризацию ручных автоморфизмов алгебры $P=\Phi[x_1, x_2]$ дает следующая
теорема.

\begin{theorem}\label{sokr}
Пусть $\phi=(f_1, f_2)$ --- ручной автоморфизм алгебры $\Phi[x_1,x_2]$. Если
$$D(\phi)\succ (x_1, x_2, 2e),$$
то автоморфизм $\phi$ является элементарно $D$-сократимым.
\end{theorem}

\begin{proof}
По следствию 1 автоморфизм
$\phi$ записывается в виде \eqref{1}.
Рассмотрим случай когда
$\lambda=id.$
Тогда
$$
\phi=\sigma_1 \circ \tau_1 \circ \sigma_2 \circ \tau_2 \circ \cdots \circ \sigma_k \circ \tau_k=(f,g).
$$
Положим
$$
\psi=\sigma_1 \circ \tau_1 \circ \sigma_2 \circ \tau_2 \circ \cdots \circ \sigma_k=(p,q).
$$
Если $\tau_k=(x_1 + h_k(x_2), x_2),$ то
$$\phi=(p+h_k(q),q)=(f,g).$$
Так как $\deg(\tau_i)=n_i$, то по предложению 1 имеем
$$
\deg(\psi)=n_1\cdots n_{k-1},\ \; \deg(\phi)=n_1\cdots n_{k-1} n_k,
$$
\begin{center}
$\deg(\psi)<\deg(\phi)$ и $D(\psi)\prec D(\phi).$\\
\end{center}
Поскольку
$$
\phi \to \psi,
$$
то автоморфизм $\phi$ является элементарно $D$-сократимым.

Допустим, что
$$
\lambda=(a_1x_1+b_1x_2+c_1, a_2x_1+b_2x_2+c_2)\neq id.
$$
Положим
$$
\omega= \sigma_1 \circ \tau_1 \circ \sigma_2 \circ \tau_2 \circ \cdots \circ \sigma_k \circ \tau_k=(p+h_k(q),q)=(r,s).
$$
Тогда по предложению 1 $\deg(r)=n_1n_2\cdots n_{k} > \deg(s)=n_1n_2\cdots n_{k-1}$.
Следовательно,
$$
\phi=\omega \circ \lambda=(a_1r+b_1s +c_1, a_2r+b_2s+c_2)=(f,g).
$$
Более того,
$$
r'=h'_k(s').
$$

Если $a_1, a_2 \neq 0,$ то $\overline{f}=\overline{g}$ и по лемме 1 автоморфизм $\phi$ элементарно $D$-сократим.

Если $a_1=0$, то $f'=b_1s',   \;g'=a_2r'=a_2 h'_k(s)$ и $b_1, a_2 \in \Phi^*.$
Тогда имеем, что
$$
D(\phi)=(\overline{s}^{n_k}, \overline{s}, |b_1 lc(g)|+|a_2 lc(f)|)=(\overline{s}^{n_k}, \overline{s}, |lc(g)|+| lc(f)|).
$$

В этом случае автоморфизм $\psi=(f, g - a_2 h_k({b_1}^{-1}f))$ является элементарным $D$-сокращением $\phi.$

Если $a_2=0,$ то $f'=a_1r'=a_1 h'_k(s), g'=b_2s'$ и $b_2, a_1 \in \Phi^*.$
Этот случай симметричен предыдущему.
\end{proof}

\begin{corollary}
Ручные и дикие автоморфизмы алгебры $\Phi[x_1, x_2]$ над конструктивным евклидовым кольцом $\Phi$ алгоритмически распознаваемы.
\end{corollary}

\begin{proof}
Проведем индукцию по показателю $D(\phi)$ автоморфизма ${\phi \in Aut(\Phi[x_1,x_2])}.$
Если $D(\phi)=(x_1, x_2, 2e),$ то $\phi$ имеет вид (2) или (3). Следовательно, автоморфизм $\phi$ является линейным. Так как любая обратимая матрица над евклидовым кольцом является произведением элементарных и диагональных матриц \cite{24}, то линейные автоморфизмы являются ручными.

Если автоморфизм $\phi$ не является элементарно $D$-сократимым, то $\phi$ является диким по теореме 2. Если $\phi$  элементарно $D$-сократим, по следствию 2 эффективно найдется $\phi'$ такой, что $D(\phi')\prec D(\phi).$ Очевидно, $\phi$ является ручным тогда и только тогда, когда $\phi'$ является ручным автоморфизмом. Так как $D(\phi')\prec D(\phi),$ то индуктивное предположение завершает доказательство.
\end{proof}
\section{Аналог автоморфизма Нагаты}
\hspace*{\parindent}
Пусть $\Phi$ --- евклидово кольцо и $0\neq z\in \Phi\setminus \Phi^*.$ Положим, что ${K=Q(\Phi)}$ --- поле частных кольца $\Phi.$
Рассмотрим следующую последовательность элементарных преобразований алгебры $K[x_1,x_2]$ над полем частных $K:$

\[(x_1,x_2)\to (zx_1,x_2)\to (zx_1-x_2^2,x_2)\to (zx_1-x_2^2,x_2+z(zx_1-x_2^2))\]

\[ \to (zx_1-x_2^2+(x_2+z(zx_1-x_2^2))^2,x_2+z(zx_1-x_2^2))=\]

\[=(z(x_1+2x_2(zx_1-x_2^2)+z(zx_1-x_2^2)^2),x_2+z(zx_1-x_2^2))\]

\[ \to (x_1+2x_2(zx_1-x_2^2)+z(zx_1-x_2^2)^2,x_2+z(zx_1-x_2^2))=\sigma. \]

Мы получили автоморфизм Нагаты $\sigma$ и он является ручным над $K$. Заметим, что $\sigma$ является автоморфизмом алгебры $\Phi[x_1, x_2]$. В своей известной работе \cite{3} М. Нагата доказал, что $\sigma$ является диким автоморфизмом алгебры $\Phi[x_1, x_2]$, если $\Phi={\bf k}[z]$ --- кольцо многочленов от одной переменной $z$. Следующее следствие является обобщением этого результата на случай произвольного евклидова кольца $\Phi$.

\begin{corollary} Автоморфизм Нагаты
$$
\sigma = (f,g)=(x_1 +2x_2w +zw^2,x_2 +zw), \:w=zx_1-x_2^2, \:0\neq z\in \Phi\setminus \Phi^*
$$
как автоморфизм алгебры $\Phi[x_1,x_2]$ над $\Phi$ является диким.
\end{corollary}
\begin{proof}
Достаточно показать, что $\sigma$ не является элементарно $D$-сократимым. Допустим, что $\sigma$ является элементарно $D$-сократимым. Тогда по следствию 2, $f'\in \left\langle g' \right\rangle$ или $g'\in \left\langle f' \right\rangle$.
Так как
$$f'=zx_2^4,\;\  g'=zx_2^2,$$
то
$f'\notin \left\langle g' \right\rangle$ и $g'\notin \left\langle f' \right\rangle$.
Следовательно, автоморфизм $\sigma$ не является элементарно $D$-сократимым и по теореме 2 не является ручным.
\end{proof}

Пусть $A=\Phi\left\langle x_1, x_2\right\rangle$ --- свободная неассоциативная алгебра с множеством свободных порождающих $X=\left\{x_1, x_2\right\}$ над евклидовым кольцом $\Phi$.
Рассмотрим автоморфизмы алгебры $A.$ Пусть $0\neq z\in \Phi\setminus \Phi^*$.
Применяя те же преобразования, что и при построении автоморфизма $\sigma$, в алгебре $A$ получаем
\begin{eqnarray}
(x_1,x_2)\to (zx_1,x_2)\to (zx_1-x_2^2,x_2)\to (zx_1-x_2^2,x_2+z(zx_1-x_2^2)) \nonumber\\
\to(zx_1-x_2^2+(x_2+z(zx_1-x_2^2))^2,x_2+z(zx_1-x_2^2))=\nonumber\\
(z(x_1+zx_2x_1-x_2x_2^2+zx_1x_2-x_2^2x_2+z^3{x_1}^2-z^2x_1x_2^2-z^2x_2^2x_1+zx_2^2x_2^2),
\nonumber \\
x_2+z(zx_1-x_2^2))\nonumber \\
\to(x_1+zx_2x_1-x_2x_2^2+zx_1x_2-x_2^2x_2+z^3{x_1}^2-z^2x_1x_2^2-z^2x_2^2x_1+zx_2^2x_2^2, \nonumber \\x_2+z(zx_1-x_2^2)) \nonumber \\
=(x_1+x_2(zx_1-x_2^2)+(zx_1-x_2^2)x_2+z(zx_1-x_2^2)^2, x_2+z(zx_1-x_2^2)). \nonumber
\end{eqnarray}
Это дает эндоморфизм
$$
\eta=(x_1+x_2(zx_1-x_2^2)+(zx_1-x_2^2)x_2+z(zx_1-x_2^2)^2,\\x_2+z(zx_1-x_2^2))
$$
алгебры $A.$

\begin{lemma} Эндоморфизм $\eta$ является автоморфизмом алгебры $A=\Phi\left\langle x_1,x_2\right\rangle$.
\end{lemma}

\begin{proof}
Пусть
$
\eta=(b_1,b_2).
$
Подалгебра $\Phi\left\langle b_1,b_2\right\rangle$, порожденная элементами $b_1,$ $b_2$ над $\Phi$, является свободной, так как линейные части элементов $b_1$, $b_2$ порождают свободную алгебру.
Докажем, что $\Phi\left\langle x_1,x_2\right\rangle = \Phi \left\langle b_1, b_2\right\rangle$. Для этого достаточно показать, что элементы $b_1, b_2$ порождают всю алгебру $\Phi\left\langle x_1, x_2\right\rangle$.
Непосредственные вычисления дают
\[zb_1-b_2^2=zx_1-x_2^2=s_1,\]
\[b_2-z(zb_1-b_2^2)=x_2=s_2,\]
\[zb_1-b_2+(b_2-z(zb_1-b_2^2))^2=zx_1=s_3,\]
\[b_1-s_2s_3+s_2s_2^2-s_3s_2+s_2^2s_2-zs_3^2+zs_3s_2^2+zs_2^2s_3-zs_2^2s_2^2=x_1=s_4.\]

Отсюда $s_1,\: s_2=x_2,\: s_3, s_4=x_1 \in\Phi\left\langle b_1, b_2\right\rangle$,
т.е. $x_1, x_2 \in \Phi\left\langle b_1, b_2\right\rangle$. Следовательно, $\Phi\left\langle x_1,x_2\right\rangle = \Phi \left\langle b_1, b_2\right\rangle$.
\end{proof}

\begin{theorem} Автоморфизм $\eta$ алгебры $A=\Phi \left\langle x_1,x_2\right\rangle$
является диким.
\end{theorem}

\begin{proof}
Пусть $\Phi[x,y]$ --- алгебра многочленов со множеством свободных порождающих $X=\left\{x, y\right\}$ над евклидовым кольцом $\Phi.$ Рассмотрим естественный гомоморфизм
$$
\tau:A\to \Phi [x,y],
$$
где $\tau(x_1)=x$ и $\tau(x_2)=y.$ Этот гомоморфизм индуцирует гомоморфизм
$$
\tau^*:Aut(A)\to Aut(\Phi [x,y]),
$$
определенный правилом
$$
\tau^*(\psi)(\tau (f))=(\tau \circ \psi)(f),
$$
где $\psi \in Aut(A), f \in A.$

Элементарному автоморфизму алгебры $A$ соответствует элементарный автоморфизм алгебры $\Phi [x,y]$.
Так как группа ручных автоморфизмов алгебры $A$ порождается элементарными автоморфизмами алгебры $A$, то гомоморфный образ ручного автоморфизма алгебры $A$ является ручным автоморфизмом алгебры $\Phi [x,y]$, т.е. $\tau^*$ индуцирует гомоморфизм групп ручных автоморфизмов
$$
\tau^*:T(A)\longmapsto T(\Phi[x,y]).
$$
Так как
$$
\tau^*(\eta)(\tau(x_1))=(\tau \circ \eta)(x_1)=\tau(\eta(x_1))=\tau(b_1)
$$
и
$$
\tau^*(\eta)(\tau(x_2))=(\tau \circ \eta)(x_2)=\tau(\eta(x_2))=\tau(b_2),
$$
то
$$
\tau^*(\eta)=(\tau(b_1),\tau(b_2))=(x+2y(zx-y^2)+z(zx-y^2)^2,y+z(zx-y^2))=\sigma
$$
есть автоморфизм Нагаты алгебры $\Phi [x,y]$ над $\Phi$. По следствию 4 автоморфизм Нагаты $\sigma$ алгебры $\Phi [x,y]$ над $\Phi$ является диким. Следовательно, $\eta$ является диким, так как его гомоморфный образ является диким.
\end{proof}

Доказательство теоремы 3 полностью проходит и для свободных неассоциативных коммутативных алгебр. Следовательно, автоморфизм
$$
\omega=(x+2y(zx-y^2)+z(zx-y^2)^2,y+z(zx-y^2)),\ 0 \neq z\in\Phi \setminus \Phi^*,
$$
свободной неассоциативно-коммутативной алгебры $B=\Phi_C\left\langle x,y \right\rangle$ над евклидовым кольцом $\Phi$ также является диким.

\end{document}